\documentclass[
  reprint,            
  aps,
  prl,
  superscriptaddress, 
  amsmath,amssymb,
  nofootinbib
]{revtex4-2}

\usepackage[T1]{fontenc}
\usepackage[utf8]{inputenc}

\usepackage{graphicx}
\usepackage{bm}
\usepackage{booktabs}
\usepackage{tabularx}
\usepackage{appendix}
\usepackage{comment}
\usepackage{hyperref}

\graphicspath{{Figures/}}

\begin{document}

\title{Quantization of Lagrangian Descriptors}

\author{Javier Jim\'enez-L\'opez}
\affiliation{Departamento de F\'isica Te\'orica, Universidad Complutense de Madrid, E-28040 Madrid, Spain}

\author{V. J. Garc\'ia-Garrido}
\affiliation{Departamento de F\'isica y Matem\'aticas, Facultad de Ciencias, Universidad de Alcal\'a, 28805 Alcal\'a de Henares, Madrid, Spain}

\date{\today}

\begin{abstract}

We formulate Lagrangian descriptors (LDs) in the path integral framework. Averaging the classical LD over fluctuations about extremal trajectories defines a quantum LD that incorporates quantum effects. Invariant manifolds, which sharply organize classical transport, become finite-width phase space structures under quantum fluctuations, and their overlap provides a geometric mechanism consistent with tunneling as fluctuation-induced delocalization of transport barriers. We demonstrate this approach for the Hamiltonian saddle, where path integral sampling reveals manifold broadening and barrier penetration. This establishes a geometric framework for studying phase space transport and tunneling beyond the classical regime, while also providing a natural route toward the application of LDs to field theory.

\end{abstract}

\maketitle

\begingroup
\renewcommand{\thefootnote}{\fnsymbol{footnote}}
\footnotetext[1]{\href{mailto:javiej09@ucm.es}{javiej09@ucm.es}}
\footnotetext[2]{\href{mailto:vjose.garcia@uah.es}{vjose.garcia@uah.es}}
\endgroup


\section{Introduction}
\label{sec:Introduction}

The concept of \emph{phase space} \cite{nolte2010tangled} is central to modern science, providing the natural arena in which dynamical phenomena unfold. Its systematic study dates back to H.~Poincar\'e's seminal work on the three-body problem \cite{Poincare1890}, whose qualitative approach to differential equations laid the foundations of modern dynamical systems theory. By introducing geometric and topological methods—most notably through invariant manifolds, periodic orbits, and recurrence—Poincar\'e revealed that the long-term evolution of trajectories is governed by intricate structures embedded in phase space. In nonlinear systems, transport is thus mediated by these invariant phase space structures \cite{wiggins2013chaotic}.
 
Understanding and characterizing these structures requires powerful analytical and computational tools. Among the most established are Poincar\'e maps, which reduce continuous flows to discrete-time dynamical systems and expose the geometry of phase space through the organization of
fixed points, islands of stability, and chaotic seas \cite{poincare1892methodes,Henon1982,wiggins2003introduction}, and Lyapunov exponents, which quantify the rate of separation of nearby trajectories and thereby provide a rigorous measure of chaos and regularity \cite{Oseledets1968,Benettin1980}. These classical diagnostics have proven indispensable across physics \cite{HOLMES1990137}, chemistry \cite{wiggins2025}, and beyond, and have been extended to the quantum domain \cite{Graham1988,Majewski2006} through, e.g., out-of-time-order correlators \cite{Larkin1969,Maldacena2016} and quantum Lyapunov exponents \cite{Rozenbaum2017}, bridging classical dynamical systems theory and quantum chaos. Nevertheless, identifying the precise geometric structures—invariant manifolds, stable and unstable foliations—that govern global transport remains a subtle and computationally demanding task, particularly in high-dimensional systems.
 
Lagrangian descriptors (LDs) have emerged as a simple yet powerful framework that addresses this challenge by revealing phase space structures directly from trajectory information \cite{mancho2013,lopesino2017,Garcia2022a}. Defined as a scalar functional evaluated along trajectories of a dynamical system, LDs are computationally efficient and have been successfully applied across disciplines \cite{madrid2009,mendoza10,wiggins2025}. Crucially, they have also demonstrated the ability to reveal the chaotic nature of trajectories in Hamiltonian systems
\cite{jimenez2026,jimenezlopez2026trajectoryarclengthrevealschaos}, complementing and extending classical diagnostics such as Poincar\'e maps and Lyapunov exponents.
 
In quantum mechanics, transport is typically analyzed through quasi-probability distributions such as the Wigner \cite{PhysRev.40.749} and Husimi \cite{husimi} functions, through semiclassical propagators \cite{doi:10.1073/pnas.14.2.178}, or Entangled Trajectory Molecular Dynamics (ETMD) \cite{PhysRevLett.87.223202, https://doi.org/10.1002/qua.10377,10.1063/1.1597496,D5CP04657B}, which have been fruitfully applied within a phase-space framework across numerous settings \cite{https://doi.org/10.1002/qute.202100016,doi:10.1142/5287}. While these approaches successfully capture interference and tunneling phenomena \cite{CohenTannoudji1977}, they do not directly identify the underlying geometric phase space structures. A geometric formulation of quantum transport barriers—analogous to the invariant manifolds of classical mechanics—thus remains largely unexplored.
 
In this Letter, we present a quantum formulation of LDs based on Feynman's path integral \cite{feynman2010quantum}. By averaging the classical LD over quantum fluctuations around the classical solution, we incorporate quantum effects that endow classical invariant manifolds with a finite width, transforming sharp transport barriers into geometric objects with intrinsic thickness. This delocalization effect yields a direct phase-space interpretation of tunneling: classically disjoint regions become connected through the overlap of broadened manifolds. We demonstrate this framework for a paradigmatic system—the Hamiltonian saddle, where the width of the quantum-delocalized invariant manifolds admits an analytical estimate controlled by the fluctuation spectrum. These results establish a geometric framework for quantum transport that bridges dynamical systems theory and path-integral quantum mechanics, and naturally suggest extensions of LDs to classical and quantum field theories.

\section{Lagrangian descriptors}
\label{sec:LDs}

The method of Lagrangian descriptors (LDs) is a trajectory-based technique from Dynamical Systems Theory \cite{Garcia2022a} that was originally developed in the field of Geophysics to analyze Lagrangian transport and mixing processes in the ocean and the atmosphere \cite{madrid2009,mendoza10}. In recent years, it has been show to reveal phase space structures-such as equilibria, stable and unstable manifolds, invariant tori, and periodic orbits-that are close related to the dynamical behaviour of trajectories in phase space \cite{mancho2013,lopesino2017}.

For a Hamiltonian system $H(\mathbf{x})$ with $n$ DoF, the evolution of the state vector $\mathbf{x} = (\mathbf{q},\mathbf{p}) \in \mathcal{S}$ is governed by Hamilton's equations of motion:
\begin{equation} \label{eq:eq_motion}
    \dot{\mathbf{x}} = \mathbf{f}(\mathbf{x}) = J \, \nabla H(\mathbf{x}) \, ,
\end{equation}
where $J$ is the Poisson matrix. A Lagrangian descriptor $\mathcal{L}$ is a non-negative scalar functional that depends on a non-negative function $\mathcal{F}(\mathbf{x}(t;\mathbf{x}_0),t)$ and the initial condition $\mathbf{x}_0$ at time $t_0$ and which is defined as:
\begin{equation} \label{eq:LD}
    \mathcal{L}(\mathbf{x_0}, t_0, T) = \int_{t_0 - T}^{t_0 + T} \text{d} t \,  \mathcal{F}(\mathbf{x}(t;\mathbf{x}_0),t)  \, .
\end{equation}
In the present work we will adopt the following integrand to define the LD functional:
\begin{equation}\label{eq:LD_func}
    \mathcal{F}(\mathbf{x},t) = \sum_{i = 1}^{2n} |f_i(\mathbf{x},t)|^{1/2} \, ,
\end{equation}
as it has been shown in the literature to be very effective for unveiling phase space structures \cite{lopesino2017}. 

\section{Theoretical formalism} 
\label{sec:theo}

Considering a Hamiltonian function of the form described in the previous section, that being:
\begin{equation} \label{eq:hamiltonian}
H(\mathbf{q},\mathbf{p}) = \frac{1}{2}\mathbf{p}^{\mathsf T} M^{-1}\mathbf{p} + V(\mathbf{q}) \, ,
\end{equation}
where $M$ is a symmetric positive-definite mass matrix and $V(\mathbf{q})$ is the potential, the corresponding configuration space action is:
\begin{equation}
S[\mathbf{q}] = \int_{-T}^{T} \text{d}t\,
\left[
\frac{1}{2}\dot{\mathbf{q}}^{\mathsf T} M \dot{\mathbf{q}} - V(\mathbf{q})
\right]\, .
\end{equation}
Letting $\mathbf{q}_{\mathrm{cl}}(t)$ be the classical trajectory satisfying the Euler-Lagrange equation, we introduce fluctuations around this trajectory as \cite{feynman2010quantum, doi:10.1142/5057,schulman2012techniques}:
\begin{equation}
\mathbf{q}(t) = \mathbf{q}_{\mathrm{cl}}(t) + \boldsymbol{\eta}(t) \, , \qquad \text{with}
\qquad \boldsymbol{\eta} \,(\pm T) = 0 \, .
\end{equation}
The action difference functional, which encodes all the physical processes of the system, is then given by:
\begin{equation}
\label{eq:deltaS}
\Delta S[\boldsymbol{\eta}; \mathbf{q}_{\mathrm{cl}}] := S[\mathbf{q}_{\mathrm{cl}} + \boldsymbol{\eta}] - S[\mathbf{q}_{\mathrm{cl}}] \, ,
\end{equation}
where $S[\mathbf{q}_{\mathrm{cl}}]$ is the classical action; and the path integral associated with the fluctuations is:
\begin{equation}
\int_{\mathcal{Y}} \mathcal{D}\boldsymbol{\eta}\;
\exp\!\left( \frac{i}{\hbar} \Delta S[\boldsymbol{\eta}; \mathbf{q}_{\mathrm{cl}}] \right) \, ,
\end{equation}
where $\mathcal{Y}$ denotes the space of real paths satisfying Dirichlet boundary conditions. Since the integrand is oscillatory, a well-defined mathematical formulation is achieved by deforming the integration contour into a complexified space. In this framework, the integral is expressed as a sum over Lefschetz thimbles \cite{Witten:2010cx, Witten:2010zr, PhysRevD.88.051501,TANIZAKI2014250,PhysRevD.111.083524, FELDBRUGGE2023169315, PhysRevD.111.085027}:
\begin{equation}
\label{eq:thimble_decomposition}
\int_{\mathcal{Y}} \mathcal{D}\boldsymbol{\eta}\; e^{\frac{i}{\hbar}\Delta S[\boldsymbol{\eta}; \mathbf{q}_{\mathrm{cl}}]}
=
\sum_{\sigma} n_\sigma
\int_{\mathcal{J}_\sigma} \mathcal{D}\boldsymbol{\eta}\; e^{\frac{i}{\hbar}\Delta S[\boldsymbol{\eta}; \mathbf{q}_{\mathrm{cl}}]} \, ,
\end{equation}
where $\mathcal{J}_\sigma$ are the steepest descent cycles associated with the relevant complex saddle points and $n_\sigma$ are intersection numbers.

Then, for any observable $\mathcal{L}[\mathbf{q},\mathbf{p}]$ defined as a scalar functional over the trajectory, the corresponding expectation value conditioned on the sector associated with $\mathbf{q}_{\mathrm{cl}}$ is defined as:
\begin{equation}
\label{eq:quantum_expec}
\langle\mathcal{L} \rangle
=
\frac{
\sum_{\sigma} n_\sigma
\int_{\mathcal{J}_\sigma} \mathcal{D}\boldsymbol{\eta}\;
\mathcal{L}[\mathbf{q}, \mathbf{p}]\;
e^{\frac{i}{\hbar}\Delta S[\boldsymbol{\eta}; \mathbf{q}_{\mathrm{cl}}]}
}{
\sum_{\sigma} n_\sigma
\int_{\mathcal{J}_\sigma} \mathcal{D}\boldsymbol{\eta}\;
e^{\frac{i}{\hbar}\Delta S[\boldsymbol{\eta}; \mathbf{q}_{\mathrm{cl}}]}
}\,.
\end{equation}

In order to analyze how the invariant manifolds of the system are broadened due to the fluctuations, we have to use a transverse coordinate. Assuming that the structure is locally characterized by the vanishing of a smooth scalar function:
\begin{equation}
    g(\mathbf{q},\mathbf{p},t) = 0\, ,
\end{equation}
one can define a transverse coordinate $u(\mathbf{q},\mathbf{p},t)$ as a function which is zero over the manifold and such that $|u(\mathbf{q},\mathbf{p},t)|$ is a perpendicular distance to the manifold. A natural way of doing this is to define the variable \cite{wiggins2003introduction}:
\begin{equation}
    u(\mathbf{q},\mathbf{p},t) = \dfrac{g(\mathbf{q},\mathbf{p},t)}{\| \nabla g(\mathbf{q},\mathbf{p},t)\|} \, ,
\end{equation}
and calculate its variance within the path integral formalism as:
\begin{equation}
    \sigma^2_u(t) 
=
\frac{
\sum_{\sigma} n_\sigma
\int_{\mathcal{J}_\sigma} \mathcal{D}\boldsymbol{\eta}\;
|u(\mathbf{q},\mathbf{p},t)|^2\;
e^{\frac{i}{\hbar}\Delta S[\boldsymbol{\eta}; \mathbf{q}_{\mathrm{cl}}]}
}{
\sum_{\sigma} n_\sigma
\int_{\mathcal{J}_\sigma} \mathcal{D}\boldsymbol{\eta}\;
e^{\frac{i}{\hbar}\Delta S[\boldsymbol{\eta}; \mathbf{q}_{\mathrm{cl}}]}
}\,,
\end{equation}
which yields the broadening of the invariant manifolds:
\begin{equation}
    \sigma_{\text{rms}} = \sqrt{\dfrac{1}{2T} \int_{-T}^{T} \text{d}t \, \sigma^2_u(t)} \, .
\end{equation}

\section{The Hamiltonian saddle with 1 DoF} \label{sec:saddle_system}

Consider an arbitrary Hamiltonian vector field with 1 DoF, and suppose it has a saddle equilibrium point at the origin at which the matrix associated with the linearization of Hamilton's equations about the equilibrium point has eigenvalues $\pm \lambda$. Using Hamiltonian normal form theory \cite{arnold1988mathematical}, one can construct the Hamiltonian function that describes the linearized dynamics of the system in phase space about the saddle equilibrium:
\begin{equation} \label{eq:saddle}
H(q,p) = \frac{\lambda}{2}\left(p^{2} - q^{2}\right) \, ,
\qquad \lambda > 0 \, .
\end{equation}
The invariant stable ($\mathcal{W}^s$) and unstable ($\mathcal{W}^u$) manifolds emanating from the saddle equilibrium act as transport barriers in phase space, and are given by the sets:
\begin{equation}
\mathcal{W}^{u/s} = \left\lbrace (q,p) \in \mathbb{R}^2 \; \Big| \; p = \pm q \right\rbrace \, 
\end{equation}
where the plus and minus signs correspond, respectively, to the unstable and stable manifolds.


The real-time phase-space action, in natural units $(\hbar = 1)$, is:
\begin{equation}
    S[q,p] = \int_{-T}^{T} \text{d}t \left(p\dot{q} - H(p,q) \right) \, ,
\end{equation}
from which the configuration-space action is recovered after integrating the momentum dependency:
\begin{equation}
    S[q,\dot{q}] = \int_{-T}^{T} \text{d}t \left( \dfrac{1}{2\lambda} \dot{q}^2 + \dfrac{\lambda}{2} q^2 \right) \, .
\end{equation}
Decomposing the path as a sum of the classical trajectory and fluctuations that satisfy Dirichlet boundary conditions, the action can be separated as:
\begin{equation}
    S[q, \dot{q}] = S[q_{\mathrm{cl}}, \dot{q}_{\mathrm{cl}}] + S[\eta, \dot{\eta}] \, ,
\end{equation}
where $S[q_{\mathrm{cl}}, \dot{q}_{\mathrm{cl}}]$ corresponds to the classical action and $S[\eta, \dot{\eta}]$ is the action of the fluctuations, given by:
\begin{equation}
    S[\eta, \dot{\eta}] = \dfrac{1}{2\lambda} \int_{-T}^{T} \text{d}t \left( \dot{\eta}^2 + \lambda^2 \eta^2 \right) \, ,
\end{equation}
which is directly linked to the following Sturm-Liuville operator:
\begin{equation} \label{eq:operator_saddle}
    \mathcal{O} = -\dfrac{\text{d}^2}{\text{d}t^2} + \lambda^{2} \, .
\end{equation}
This allows us to express the fluctuations in terms of the eigenfunctions of $\mathcal{O}$ as:
\begin{equation} \label{eq:spectral_expasion}
    \eta(t) = \sum_{n = 1}^{\infty} c_n \phi_n(t) \, , 
\end{equation}
where $c_n$ are the Fourier coefficients and the eigenfunctions of $\mathcal{O}$ are:
\begin{equation}
    \phi_{n}(t) = \dfrac{1}{\sqrt{T}} \sin (k_n (t + T)) \, ,\qquad k_n = \dfrac{n\pi}{2T} \, .
\end{equation}
Since the eigenfunctions are orthogonal, the action of the fluctuations can be rewritten as:
\begin{equation}
    S[\eta, \dot{\eta}] = \dfrac{1}{2\lambda} \sum_{n=1}^{\infty} (k^{2}_{n} + \lambda^2) c^2_n \, .
\end{equation}

To make the integral converge, we deform each mode into its steepest-descent contour by letting \cite{FELDBRUGGE2023169315}:
\begin{equation}
    c_n = e^{i\pi / 4} y_n \, \qquad \text{with} \qquad y_n \in \mathbb{R} \, ,
\end{equation}
so that the oscillatory weight becomes a real Gaussian:
\begin{equation}
    \exp \left[ \dfrac{i}{2\lambda} \left( k^2_n + \lambda^2 \right) c^2_n \right] \rightarrow \exp \left[ - \dfrac{1}{2\lambda} \left( k^2_n + \lambda^2 \right) y^2_n \right] \, ,
\end{equation}
and the fluctuation field becomes complex:
\begin{equation}
    \eta(t) = e^{i\pi / 4} \sum_{n = 1}^{\infty} y_n \phi_n(\tau) \, .
\end{equation}
This method allows for the dynamical evolution to take place in real time while the integration contour lives in a complexified space. Then, the thimble variables $y_n$ are independent Gaussians with zero mean and variance:
\begin{equation}
     \text{Var}(y_n) = \sigma_n^2 = \dfrac{\lambda}{k^2_n + \lambda^2} \, .
\end{equation}

Since we are interested in estimating the width that the invariant manifolds acquire due to the quantum fluctuations, and taking into account the symmetry of the system, we define the new coordinate:
\begin{equation}
    u(t) = \dfrac{p(t) - q(t)}{\sqrt{2}} \, ,
\end{equation}
so that the fluctuation in the transverse direction is:
\begin{equation}
        \delta u(t) = \dfrac{1}{\sqrt{2}} \left( \dfrac{1}{\lambda} \dot{\eta}(t) - \eta(t) \right) = \dfrac{e^{i\pi / 4}}{\sqrt{2}} \sum_{n=1}^{N} y_n a_n(t) \, ,
\end{equation}
where we have replaced the infinite sum for a finite one to avoid divergences and with:
\begin{equation}
    a_n(t) = \dfrac{1}{\lambda} \dot{\phi}_n(t) - \phi_n(t) \, .
\end{equation}

\begin{figure*}[htpb]
    \centering
    A) \includegraphics[scale = 0.5]{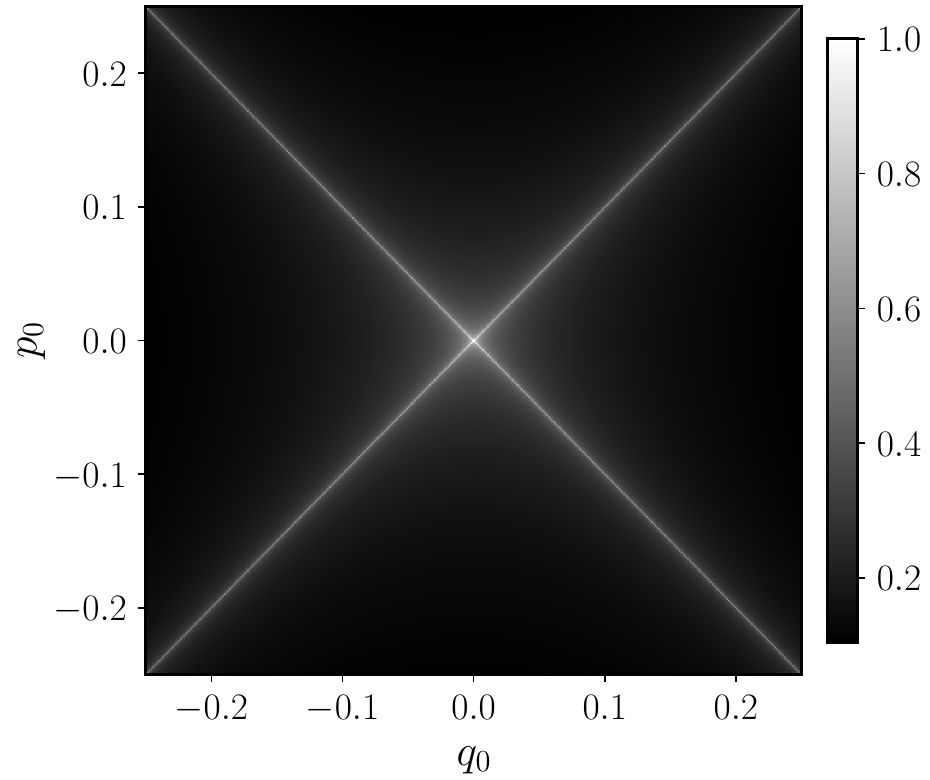}
    B) \includegraphics[scale = 0.5]{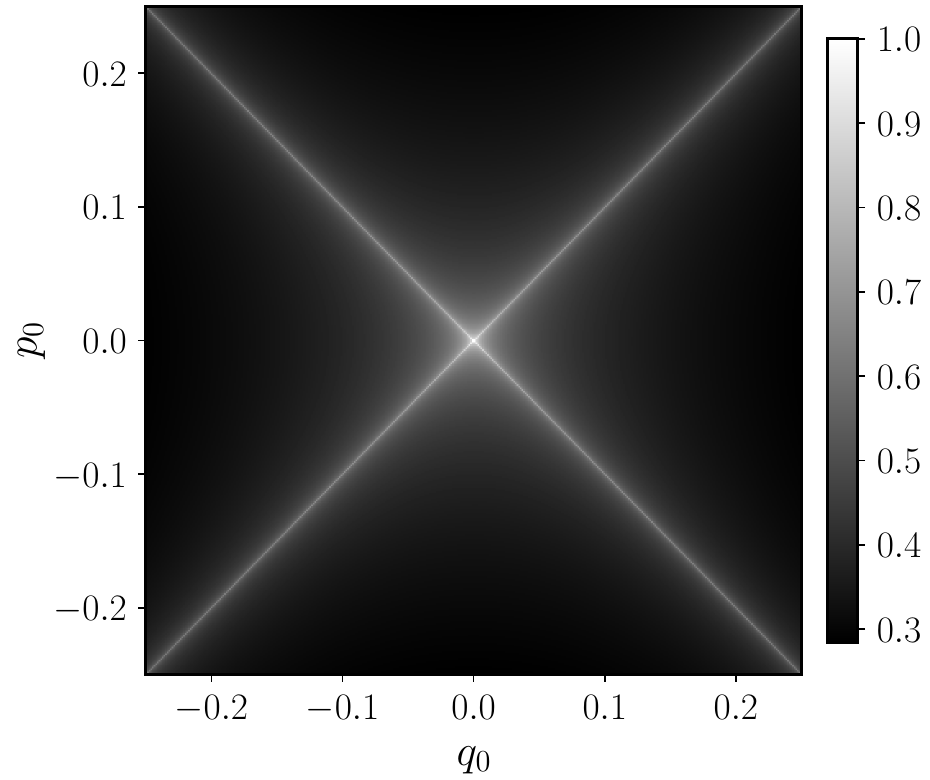}
    \caption{Difference between the classical and quantum Lagrangian descriptor for the Hamiltonian saddle given in Eq.~\eqref{eq:saddle} calculated with A) $10$ modes and B) $800$ modes in the decomposition of the operator in Eq.~\eqref{eq:operator_saddle}. As predicted by Eq.~\eqref{eq:width}, the broadening of the manifolds grows as the number of modes considered in the calculation becomes larger. The numerical calculation was done with $\lambda = 3$, $T = 8/3$ and the initial conditions where sampled from a $400\times 400$ regular grid with $1200$ samples per initial condition. Note that while the quantum LD faithfully encodes the effective broadening through the path-integral average, the visual contrast is constrained by the finite grid resolution and colormap dynamic range; structures narrower than the grid spacing are captured analytically by Eq.\eqref{eq:width} but cannot be visually resolved.}
    \label{fig:saddle}
\end{figure*}

Since the coefficients $y_n$ are independent, the variance of the fluctuation in the transverse direction is:
\begin{equation}
    \langle |\delta u(t)|^2 \rangle = \dfrac{1}{2} \,\sum_{n = 1}^{N} \sigma_n^2 a^2_n(t) = \dfrac{\lambda}{2} \sum_{n= 1}^{N}\dfrac{a^2_n(t)}{k^2_n + \lambda^2} \, ,
\end{equation}
and then, the width of the invariant manifolds can be estimated as:
\begin{equation} \label{eq:width}
    \sigma_{\text{rms}} = \sqrt{\dfrac{1}{2T} \int_{-T}^{T} \text{d}t \, \langle |\delta u(t)|^2 \rangle} = \sqrt{\dfrac{N}{4T\lambda}} \, .
\end{equation}
Due to the symmetry of the invariant manifolds in this system, the same result holds for the other manifold.

The expression in Eq.~\eqref{eq:width} implies that the broadening of the invariant manifolds grows with the number of modes considered in the expansion, reflecting the well-known ultra-violet sensitivity of path integrals when observables depend on the derivatives of the path. This behaviour is not pathological but rather a standard feature of functional integrals, whose continuum definition requires regularization \cite{Bastianelli_van_Nieuwenhuizen_2006, PhysRevD.106.126002}. In practice, the path integral is always defined through a discretization of the time interval where the time-step size provides a natural UV cut-off. The truncation of the spectral expansion performed in Eq.~\eqref{eq:spectral_expasion} plays the role of a regulator \cite{PhysRevD.58.044002} analogous to the lattice spacing in lattice field theory \cite{Wilson:2004de} or the momentum cutoff commonly used in QFT \cite{Peskin:1995ev}. Observables computed from the path integral are interpreted at finite resolution, as in lattice MC implementations of quantum mechanics \cite{Gattringer:2010zz,frenkel2023understanding} and statistical field theory \cite{binder1992monte, newman1999monte}. Such regularizations are ubiquitous in the path integral formalism and provide well-defined numerical and analytical results. In this sense, the finite width obtained for the invariant manifolds should be understood as a regulator dependent effective broadening induced by the quantum fluctuations at the chosen resolution scale. Physical predictions can be extracted by noting that the \emph{ratio} of broadenings between two systems with parameters $(\lambda_1,T_1)$ and $(\lambda_2,T_2)$ evaluated at the same $N$,
\begin{equation}\label{eq:ratio}
  \frac{\sigma_\mathrm{rms}^{(1)}}{\sigma_\mathrm{rms}^{(2)}} = \sqrt{\frac{\lambda_2\, T_2}{\lambda_1\, T_1}}\,,
\end{equation}
is exactly cut-off independent and constitutes a genuine physical prediction, allowing the broadening of a family of systems to be analyzed based on a reference one. 


In Fig.~\ref{fig:saddle} we display the normalized difference between the quantum LD, obtained by averaging the functional given in Eq.~\eqref{eq:LD}, where the modulus is defined as $|a| = \sqrt{a a^*}$ in order to guarantee that the LD is a real and non-negative scalar functional, with the input function of Eq.~\eqref{eq:LD_func}, and the classical LD, calculated as:
\begin{equation}
    \mathcal{L} = \int_{-T}^{T} \,  \text{d}t \, \left( |q_{\mathrm{cl}}(t)|^{1/2} +  |p_{\mathrm{cl}}(t)|^{1/2}\right) \, .
\end{equation}

\begin{figure}[htbp]
    \centering
    \includegraphics[scale = 0.41]{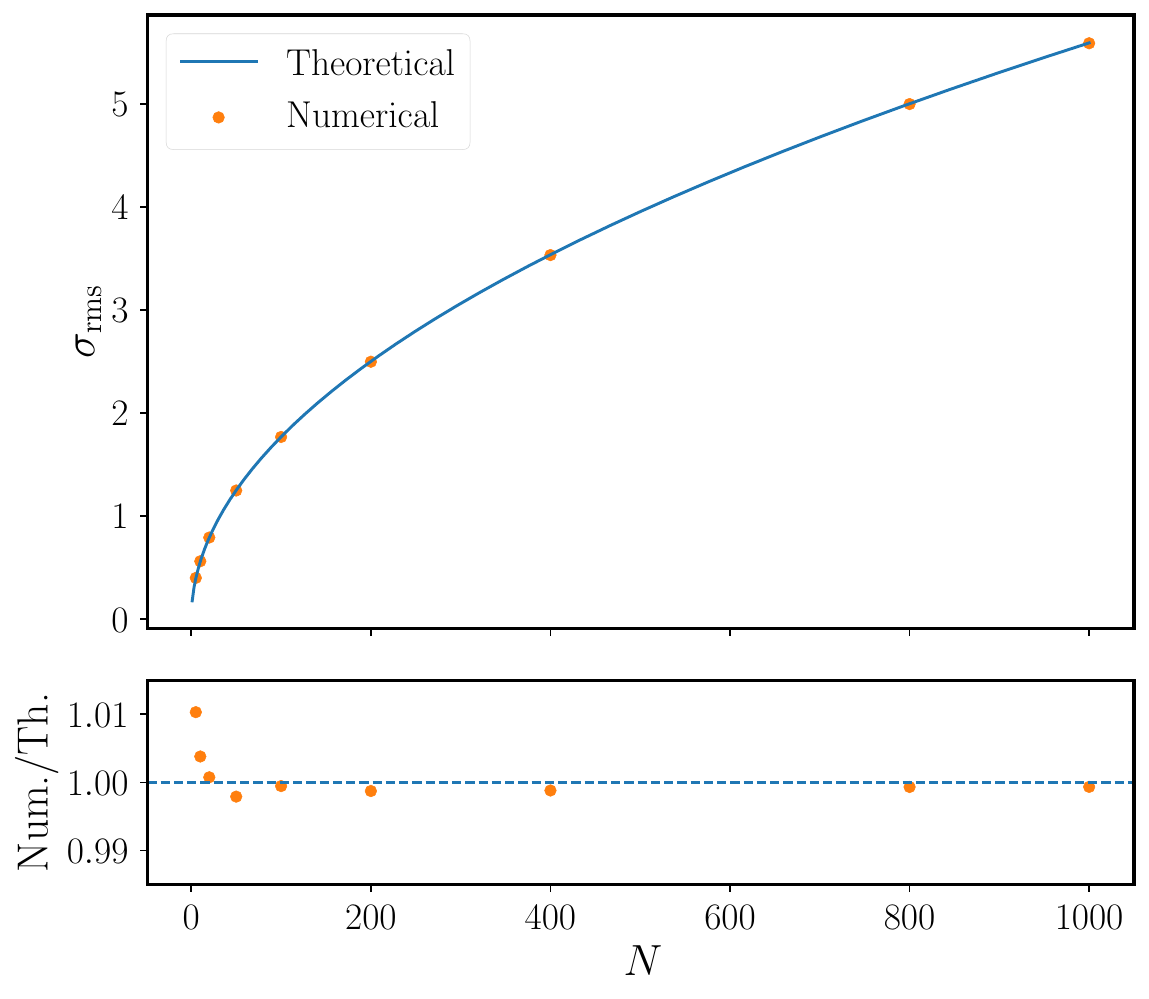}
    \caption{Width of the invariant manifolds computed as a function of the number of modes $N$ for $1200$ samples and $T\lambda = 8$. In blue we depict the theoretical curve given in Eq.~\eqref{eq:width} and in orange the result from the simulation. The agreement between the theoretical prediction and the numerical result is within $1\%$.}
    \label{fig:sigma}
\end{figure}

As predicted by Eq.~\eqref{eq:width}, the invariant manifolds acquire a broadened structure that grows the more modes are considered in the spectral decomposition of the operator given in Eq.~\eqref{eq:operator_saddle}. This effect can be visually appreciated by comparing  Fig.~\ref{fig:saddle} A), obtained with $10$ modes, and Fig.~\ref{fig:saddle} B), calculated with $800$ modes. However, while the quantum LD faithfully encodes the effective broadening of the invariant manifolds through the path-integral average, the visual contrast in the figures is inherently limited by the finite sampling resolution: structures whose width falls below the grid spacing are not visually resolved even though they are taken into account in the system's evolution.

The quantitative agreement between the analytical prediction of Eq.~\eqref{eq:width} and the Monte Carlo estimate of $\sigma_{\mathrm{rms}}$ is shown in Fig.~\ref{fig:sigma} as a function of the number of modes $N$, for fixed $T\lambda = 8$ and $1200$ samples per initial condition, resulting in an agreement within $1\%$ over the entire range explored.

\section{Conclusions}
\label{sec:Conclusions}


We have introduced a quantum formulation of Lagrangian descriptors within the path integral framework, extending the geometric description of classical transport into the quantum regime. More broadly, the present work positions LDs as an observable that admits a meaningful formulation in both classical and quantum dynamics. In this setting, classical invariant manifolds are broadened by quantum fluctuations, acquiring a non-zero width that can be characterized analytically for a Hamiltonian saddle. This broadening leads to the overlap of structures that are classically separated, consistent with tunneling as fluctuation-induced delocalization of invariant manifolds. These results open a route toward geometric characterization of quantum transport phenomena and establish new, concrete connections between quantum mechanics and the theory of dynamical systems. 

Furthermore, because the present construction is formulated in the path-integral formalism, it naturally suggests extensions of LDs to classical and quantum field theories. Identifying the associated invariant structures and computing the corresponding manifold broadening in field-theory phase space would generalize the geometric description of transport barriers to an infinite-dimensional setting, opening new directions in areas such as quantum field theory, statistical field theory, and cosmology.

\section{Acknowledgments}

The authors acknowledge Prof.~Jos\'{e} Bienvenido S\'{a}ez Landete (Universidad de Alcal\'{a}) for providing the computational resources used in this work. Javier Jim\'{e}nez-L\'{o}pez also acknowledges the valuable discussions with Prof.~Miguel \'{A}ngel Hidalgo Moreno (Universidad de Alcal\'{a}), which contributed to improving the results presented in this paper.

\appendix

\bibliography{referencias}

\end{document}